\documentclass[letterpaper, 10 pt, conference]{ieeeconf}  

\IEEEoverridecommandlockouts                              

\usepackage{amsmath,amssymb}
\usepackage{graphics,graphicx,epsfig,algorithm,algpseudocode}
\usepackage{subcaption}
\usepackage{booktabs}

\usepackage{acronym}
\usepackage{bookmark}
\usepackage{adjustbox}
\usepackage{csquotes}
\usepackage{dsfont}
\usepackage[font=small,labelfont=bf]{caption}
\usepackage{todonotes}
\usepackage{mathtools}
\usepackage{tcolorbox}
\usepackage{hyperref}

\newtheorem{definition}{Definition}
\newtheorem{theorem}{Theorem}

\newtheorem{lemma}{Lemma}

\newtheorem{assumption}{Assumption}

\title{\LARGE \bf
Hyperexponential Stabilization of Nonholonomic Systems\\ by Oscillating Controls
}

\author{Victoria Grushkovskaya$^{1,3}$, Mark\'eta Hrani\v{c}kov\'{a}$^{1}$,  and Alexander Zuyev$^{2,3}$
\thanks{$^{1}$Department of Mathematics, University of Klagenfurt, 9020 Klagenfurt am W\"orthersee, Austria
        {\tt\footnotesize \{viktoriia.grushkovska,marketa.hranickova\}@aau.at}}%
    \thanks{$^{2}$Max Planck Institute for Dynamics of Complex Technical Systems, 39106 Magdeburg, Germany
       {\tt\footnotesize zuyev@mpi-magdeburg.mpg.de} }
\thanks{$^{3}$Institute of Applied Mathematics \& Mechanics, National Academy of Sciences of Ukraine}
}

\begin{document}

\maketitle
\thispagestyle{empty}
\pagestyle{empty}

\begin{abstract}
We address the problem of designing time-varying feedback controls that stabilize a given equilibrium of a nonholonomic system at a faster than the exponential rate. We consider systems with an arbitrary degree of nonholonomy and prove a theorem on practical hyperexponential stabilization under a suitable assumption on the parameterized dither signals. 
The convergence proof is based on a comparison equation with a prescribed power of logarithmic terms.
We apply this result to stabilize a model of a car with two off-hooked trailers.
\end{abstract}

\section{Introduction}
We consider driftless control-affine systems of the form
\begin{equation}
\dot x = \sum_{k=1}^m u_k f_k(x), \quad x\in D\subset \mathbb R^n,\;u\in \mathbb R^m,
\label{Sigma}
\end{equation}
where $x=(x_1,\ldots,x_n)^\top$ is the state, $u=(u_1,\ldots,u_m)^\top$ is the control, $m<n$, 
$D$ is a domain such that $0\in D$, and $f_k:D\to\mathbb R^n$ are smooth vector fields.

For such systems, the control vector fields do not span the whole tangent space, while the missing directions are generated by their Lie brackets. This makes stabilization of nonholonomic systems essentially different from that of fully actuated systems and motivates the use of time-varying and oscillating controls~\cite{bloch2015nonholonomic,coron2007control,grushkovskaya2018obstacle,grushkovskaya2024design}. In~\cite{grushkovskaya2024design}, oscillating feedback controls were constructed for exponential stabilization of driftless nonholonomic systems with an arbitrary degree of nonholonomy.

In this paper, we consider the same class of systems and investigate whether faster convergence can be achieved. For this purpose, we replace the constant control gain used in~\cite{grushkovskaya2024design} by a state-dependent gain and show that the corresponding $\pi_\varepsilon$-solutions exhibit hyperexponential decay until entering an $\varepsilon$-dependent neighborhood of the origin, whose radius tends to zero as $\varepsilon\to0$.

Hyperexponential stability and stabilization have recently received increasing attention. Faster-than-exponential and bounded time-varying stabilization of linear systems was studied in~\cite{inoue2011nonexponential,chu2022hyperexponential}, while Lyapunov characterizations and sampled-time realizations were developed in~\cite{zimenko2023hyperexponential,zimenko2025hyperexponential}. Constructive designs were obtained for integrator systems in continuous and discrete time~\cite{efimov2023doubleintegrator,wang2024hyperexponential}, including unmatched perturbations~\cite{labbadi2024doubleintegrator}, and for several classes of nonlinear and time-delay systems~\cite{efimov2024damping,labbadi2025nonlinear,nekhoroshikh2023hyperexponential,iftime2026hyperexponential}. Hyperexponential convergence has also been considered for high-order differentiators~\cite{wang2025differentiator}.

Related accelerated stabilization problems for nonholonomic systems have mainly been studied for systems in chained form. Prescribed-time stabilization of high-order and perturbed chained systems was considered in~\cite{zhang2022prescribed,zhou2023prescribed} using state transformations, time-varying high-gain feedback, and backstepping. More recently, hyperexponential stabilization of perturbed unicycle dynamics was addressed in~\cite{labbadi2025unicycle} through state and time transformations leading to interconnected integrator-type subsystems. In contrast, the approach considered here applies directly to general driftless nonholonomic systems and uses oscillating controls associated with Lie brackets of arbitrary degree, for systems which are not transformed to chained form.

The main contributions of this paper are as follows. First, we introduce scalar comparison equations with hyperexponentially decaying solutions and use them to construct a state-dependent logarithmic control gain. The proposed gain yields several types of hyperexponential decay rate. Second, we develop a constructive oscillating-feedback design for general driftless nonholonomic systems with an arbitrary degree of nonholonomy and establish hyperexponential decay of the corresponding solutions to an arbitrary small neighborhood of the origin. Third, we apply the proposed control scheme to a car with two off-hooked trailers. Stabilizing controls for this system were previously constructed in~\cite{ZG24} using a nilpotent approximation. Here, we apply the oscillating feedback to the original nonlinear system and show that it provides hyperexponential stabilization in the above sense.

\section{Preliminaries}

\subsection{Notations and Definitions}

We denote by $B_{\delta}(x)$ and $\overline{B_\delta(x)}$ the open and closed $\delta$-neighborhoods of $x\in\mathbb R^n$ with $\delta>0$, respectively.
The Lie bracket of the vector fields $f_i$ and $f_j$ is defined by $[f_i,f_j](x)=\frac{\partial f_j(x)}{\partial x}f_i(x)-\frac{\partial f_i(x)}{\partial x}f_j(x)$.
For a multi-index $I_\ell=(i_1,\ldots,i_\ell)$ with
$i_j\in\{1,\ldots,m\}$, let
$f_{I_1}=f_{i_1}$ and, for $\ell\ge2$, define
the $\ell$-fold iterated Lie bracket by
$f_{I_\ell}(x)=\left[f_{i_1},\left[f_{i_2},\left[\dots\left[f_{i_{\ell-1}},f_{i_\ell}\right]\right]\right]\right](x)$; 
for $I_1\in\{1,\ldots,m\}$, $f_{I_1}$ is a vector field in~\eqref{Sigma}.
\begin{definition}
Given a time-varying feedback $u=h(t,x)$, $h:{\mathbb R}_{\ge 0}\times D\to\mathbb R^m$, an initial condition $x^0\in D$, and  $\varepsilon>0$, a {\em $\pi_\varepsilon$-solution} of system~\eqref{Sigma} corresponding to $x^0$ and $h$ is a continuous function $x:{\mathbb R}_{\ge 0}\to D$ satisfying $x(0)=x^0$ and
$$\dot x(t)=\sum_{k=1}^m h_k(t,x(t_j))\,f_k(x(t)),
\,
t\in(t_j,t_{j+1}),$$
for each $j=0,1,2,\ldots$, with $t_j=j\varepsilon$.
\end{definition}

\subsection{Controllability Assumption}

Throughout this paper, we assume that system~\eqref{Sigma} has degree of nonholonomy $\nu\ge 2$ at every $x\in D$. We formalize this property using left-iterated Lie brackets as follows. 

\begin{assumption}\label{rank_nu}
There exist nonempty sets of multi-indices $S_\ell\subseteq\{1,\ldots,m\}^\ell$ for
$\ell=1,...,\nu$  such that $\sum\limits_{\ell=1}^\nu |S_\ell|=n$, and the $n\times n$ matrix
$
\mathcal F(x)=\Big(f_{I_\ell}(x)\Big)_{I_\ell\in S_\ell,\,
\ell=1,...,n}
$,
whose columns are the vectors  $f_{I_\ell}(x)$, is nonsingular for all  $x\in D$.
\end{assumption}
\subsection{Control Functions}
Following the approach of~\cite{grushkovskaya2024design}, we consider oscillating control functions of the form
\begin{equation}\label{cont_nu}
u^\varepsilon(t,x)=\sum_{\ell=1}^\nu\varepsilon^{\frac{1}{\ell}-1}\sum_{I_\ell\in S_\ell}\phi^{(\varepsilon)}_{I_\ell}(t,a_{I_\ell}(x)).
\end{equation}
Here, $\varepsilon>0$ is a parameter, $\phi^{(\varepsilon)}_{I_\ell}:{\mathbb R}_{\ge 0}\times{\mathbb R}\to{\mathbb R}^m$ are parameterized dither signals, and $a_{I_\ell}:D\to\mathbb R$ are state-dependent coefficients. The components of $\phi^{(\varepsilon)}_{I_\ell}$ are denoted by $\phi^{(\varepsilon,k)}_{I_\ell}$, i.e., for any $\xi\in\mathbb R$,
$\phi^{(\varepsilon)}_{I_\ell}(t,\xi)=\left(\phi^{(\varepsilon,1)}_{I_\ell}(t,\xi),\ldots,\phi^{(\varepsilon,m)}_{I_\ell}(t,\xi)\right)^\top\in\mathbb R^m$.
As in~\cite{grushkovskaya2024design}, we make the following assumption on these dither signals.

\begin{assumption}\label{as_dither}
There exists a family of parameterized dither signals $\phi^{(\varepsilon)}_{I_\ell}(t,\xi)$, $I_\ell\in S_\ell$, $\varepsilon>0$, $\xi\in\mathbb R$, satisfying the following properties for all initial data $x^0\in K\Subset D$, where $K$ is a compact containing $\overline{B_r(0)}$ for some $r>0$:
\begin{itemize}
\item[P1)]
$\phi^{(\varepsilon)}_{I_\ell}(t,\xi)$ are continuous in $(t,\xi)\in{\mathbb R}_{\ge 0}\times{\mathbb R}$, $\varepsilon$-periodic in $t$ for each fixed $\xi$, and there exist constants $C_{I_\ell}>0$ such that
$
\|\phi^{(\varepsilon)}_{I_\ell}(t,\xi)\|\le C_{I_\ell}|\xi|^{1/\ell}
$
for all $t\ge0$, $\xi\in\mathbb R$;

\item[P2)]
there exists $\hat\delta>0$ such that, for any column vector
$
\bar\xi=(\xi_{I_\ell})^T_{I_\ell\in S_\ell,\ell=1,...,\nu}\in B_{\hat\delta}(0)\subset\mathbb R^n,
$
the solution $x(t;x^0,u(\cdot))$ of system~\eqref{Sigma} with $x^0\in K$ and the control \begin{equation}
u(t)=\sum_{\ell=1}^\nu\varepsilon^{\frac{1}{\ell}-1}
\sum_{I_\ell\in S_\ell}
\phi^{(\varepsilon)}_{I_\ell}(t,\xi_{I_\ell})
\label{open_loop}
\end{equation}
admits the representation
\begin{align}\label{volt_nu}
x(\varepsilon;x^0,u(\cdot))
=
x^0+\varepsilon{\mathcal F}(x^0)\bar\xi
+\Omega(\varepsilon,x^0,\bar\xi)
\end{align}
for all sufficiently small $\varepsilon>0$;

\item[P3)]
there exist $\varpi>0$, $\eta>0$, and $\varepsilon_0>0$ such that
$
\|\Omega(\varepsilon,x^0,\bar\xi)\|
\le
\varpi\varepsilon^{1+\eta}\|\bar\xi\|
$
for all $x^0\in K$, $\bar\xi\in B_{\hat\delta}(0)$, and $\varepsilon\in(0,\varepsilon_0]$.
\end{itemize}
\end{assumption}

In~\cite{grushkovskaya2024design}, it was proved that, under the assumptions stated above, the time-varying feedback control~\eqref{cont_nu} with
$$
a(x)=\left(a_{I_\ell}(x)\right)^\top_{I_\ell\in S_\ell,\ell=1,...,\nu}
=-\gamma\mathcal F^{-1}(x)x,\; \gamma>0,
$$
exponentially stabilizes the origin of system~\eqref{Sigma}. In particular, the corresponding $\pi_\varepsilon$-solutions  satisfy
$$
\|x(t)\|=O\left(e^{-\bar\gamma t/\nu}\right)
\quad\text{as }t\to\infty.
$$
In this paper, we replace the constant gain $\gamma$ with a state-dependent gain $\gamma(x)$ and consider the control~\eqref{cont_nu} with
\begin{equation}
\label{ax}
a(x)=-\gamma(x)\mathcal F^{-1}(x)x.
\end{equation}
We establish conditions on $\gamma(x)$ under which the corresponding $\pi_\varepsilon$-solutions exhibit hyperexponential decay until they enter an $\varepsilon$-dependent neighborhood of the origin, whose radius tends to zero as $\varepsilon\to0$.

\section{Main Results}

\subsection{Comparison Equation}
\label{sec_comp}

In~\cite{grushkovskaya2024design}, exponential convergence was obtained for the class of nonholonomic systems considered in this paper. Here, our objective is to extend these results to hyperexponential convergence. To describe the corresponding convergence rates, consider the scalar comparison equation
\begin{equation}
\dot v=-\lambda v\left|\ln\frac{v}{r}\right|^\alpha,\, 0<v<r,
\label{comparison_log}
\end{equation}
where $\lambda>0$, $r>0$, and $\alpha\in(0,1]$. The right-hand side of~\eqref{comparison_log} is extended by continuity to $0$ at $v=0$.

We say that a function $v(t)>0$\textit{ converges hyperexponentially} to zero if
\begin{equation}
\lim_{t\to\infty}\frac{\ln v(t)}{t}=-\infty.
\label{hyperexp_def}
\end{equation}

The convergence properties of~\eqref{comparison_log} are summarized in the following lemma.

\begin{lemma}
\label{lem_comparison}
\emph{Let $v(0)=v_0\in(0,r)$ and $\alpha\in(0,1]$. Then the solution of~\eqref{comparison_log} converges hyperexponentially to zero. More precisely,
\begin{itemize}
\item[(a)] If $\alpha\in(0,1)$, then
$
v(t)=e^{-c t^p+o(t^p)}$ as $t\to\infty$,
where $p=\frac{1}{1-\alpha}>1$ and $c=(\lambda(1-\alpha))^{\frac{1}{1-\alpha}}$.
\item[(b)] If $\alpha=1$, then
$
v(t)=r e^{-e^{\lambda t+O(1)}}$ as $t\to\infty$.
\end{itemize}}
\end{lemma}

\begin{proof}
Let
$
y(t)=\left|\ln\frac{v(t)}{r}\right|.
$
Since $v(t)\in(0,r)$, one has $y(t)=-\ln(v(t)/r)=\ln(r/v(t))>0$, and therefore
$$
\dot y=-{\dot v}/{v}=\lambda y^\alpha,
$$
with $y(0)=y_0:=\ln(r/v_0)>0$.

If $\alpha\in(0,1)$, direct integration gives
$$
y(t)=\left(y_0^{1-\alpha}+\lambda(1-\alpha)t\right)^{\frac{1}{1-\alpha}},
$$
and therefore
$
v(t)=r e^{-\left(\left(\ln\frac{r}{v_0}\right)^{1-\alpha}+\lambda(1-\alpha)t\right)^{\frac{1}{1-\alpha}}}.
$\\
Let $p=\frac{1}{1-\alpha}>1$ and $c=(\lambda(1-\alpha))^{\frac{1}{1-\alpha}}$. Then
$$
\left(y_0^{1-\alpha}+\lambda(1-\alpha)t\right)^{\frac{1}{1-\alpha}}=ct^p+o(t^p)
$$
as $t\to\infty$. Hence,
$
v(t)=e^{-ct^p+o(t^p)}.
$
Since $p>1$, it follows that
$
\frac{\ln v(t)}{t}=-ct^{p-1}+o(t^{p-1})\to-\infty,
$
which proves hyperexponential convergence.

If $\alpha=1$, then $\dot y=\lambda y$, and hence
$
y(t)=y_0e^{\lambda t}.
$
Consequently,
$
v(t)=r e^{-y_0e^{\lambda t}}=r e^{-e^{\lambda t+\ln y_0}}=r e^{-e^{\lambda t+O(1)}}.
$
Moreover,
$
\frac{\ln v(t)}{t}=\frac{\ln r}{t}-\frac{y_0e^{\lambda t}}{t}\to-\infty\text{ as }t\to\infty,
$
so the convergence is hyperexponential also in this case.
\end{proof}

\subsection{State-Dependent Control Gain}

To construct the state-dependent gain in~\eqref{ax}, we first consider the reference system
\begin{equation}
\dot{\bar x}=-\gamma(\bar x)\bar x,\quad \bar x\in B_r(0)\subset\mathbb R^n.
\label{sys_ref}
\end{equation}
We choose $\gamma(\bar x)$ so that the trajectories of~\eqref{sys_ref} exhibit the hyperexponential convergence behavior described by the comparison equation~\eqref{comparison_log}.

Let $v=\|\bar x\|$. For $\bar x\neq0$, differentiating $v^2=\|\bar x\|^2$ along the trajectories of~\eqref{sys_ref} gives
$
\dot v=-\gamma(\bar x)v.
$
Therefore, choosing
\begin{equation}
\gamma(\bar x)=\lambda\left|\ln\frac{\|\bar x\|}{r}\right|^\alpha,\, 0<\|\bar x\|<r,
\label{gamma}
\end{equation}
where $\lambda>0$ and $\alpha\in(0,1]$, yields
$
\dot v=-\lambda v\left|\ln\frac{v}{r}\right|^\alpha.
$
Thus, by Lemma~\ref{lem_comparison}, $v(t)=\|\bar x(t)\|$ converges hyperexponentially to zero, and therefore the trajectories of the reference system~\eqref{sys_ref} converge hyperexponentially to the origin.

We now return to the oscillating feedback~\eqref{cont_nu} introduced in Section~2:
$$
u^\varepsilon(t,x)=\sum_{\ell=1}^\nu\varepsilon^{\frac{1}{\ell}-1}\sum_{I_\ell\in S_\ell}\phi^{(\varepsilon)}_{I_\ell}(t,a_{I_\ell}(x)),
$$
where the functions $\phi^{(\varepsilon)}_{I_\ell}$ satisfy properties {\rm P1)}--{\rm P3)}. The state-dependent coefficients are chosen according to~\eqref{ax},
$$
a(x)=-\gamma(x)\mathcal F^{-1}(x)x,
$$
where
$
\gamma(x)=\lambda\left|\ln\frac{\|x\|}{r}\right|^\alpha,\, 0<\|x\|<r,
$
and $a(0)=0$. We assume that $r$ and $\lambda$ are chosen so that
\begin{equation}
a(B_r(0))\subset B_{\hat\delta}(0).
\label{a_admissible}
\end{equation}

With this choice,
$
\mathcal F(x)a(x)=-\gamma(x)x,
$
so that the leading term in the representation~\eqref{volt_nu} coincides with the vector field of the reference system~\eqref{sys_ref}. Hence, the feedback~\eqref{cont_nu} is constructed to approximate the hyperexponential convergence behavior of the reference system~\eqref{sys_ref}.

\subsection{Convergence Analysis}

We now state the main stabilization result for system~\eqref{Sigma} with the control laws defined by~\eqref{cont_nu},~\eqref{ax}, and~\eqref{gamma}.

\begin{theorem}
\label{thm_main}
\textit{Let $\nu\ge 2$ and $f_k\in C^{\nu+1}(D;\mathbb R^n)$ for $k=1,\ldots,m$, and suppose that Assumption~\ref{rank_nu} holds.
Let the dither signals $\phi^{(\varepsilon)}_{I_\ell}$ satisfy Assumption~\ref{as_dither} with $K=\overline{B_r(0)}\subset D$ for some $r>0$.
Consider the time-varying feedback defined by~\eqref{cont_nu}, with $a(x)$ given by~\eqref{ax} and $\gamma(x)$ given by~\eqref{gamma}, assuming $a(0)=0$. 
Then, for every $\theta\in(0,1)$, there is an $\hat\varepsilon>0$ such that, for every $\varepsilon\in(0,\hat\varepsilon)$ and every $x_0\in B_r(0)$, there exist positive constants $T_\varepsilon$, $\bar\delta_\varepsilon$, and $C_\varepsilon$ such that the corresponding $\pi_\varepsilon$-solution of~\eqref{Sigma},~\eqref{cont_nu} satisfies:
\begin{enumerate}
\item if $\alpha\in(0,1)$, then
$\|x(t)\|\le C_\varepsilon e^{-\sigma t^{\frac{1}{1-\alpha}}}$ for $t\in[0,T_\varepsilon)$;
\item if $\alpha=1$, then
$
\|x(t)\|\le C_\varepsilon e^{-\sigma e^{\bar\lambda t}}
$
for $t\in[0,T_\varepsilon)$, with $
\bar\lambda=\frac{(1-\theta)\lambda}{2}.
$
\end{enumerate}
Furthermore,
$
\|x(t)\|\le\bar\delta_\varepsilon$  for all  $t\ge T_\varepsilon$,
where $\bar\delta_\varepsilon\to0$ as $\varepsilon\to0$.}
\end{theorem}
The explicit estimates for the decay rate and parameters  $\hat\varepsilon$, $T_\varepsilon$, and $\bar\delta_\varepsilon$ are obtained in the below proof.

\begin{proof}
To analyze the convergence of the $\pi_\varepsilon$-solutions of~\eqref{Sigma},~\eqref{cont_nu}, we consider the representation~\eqref{volt_nu}. Let $t_j=j\varepsilon$ and $x_j=x(t_j)$. By property {\rm P2)}, the solution over one sampling interval satisfies
\begin{equation}
x_{j+1}=x_j+\varepsilon\mathcal F(x_j)a(x_j)+\Omega_j
=x_j-\varepsilon\gamma(x_j)x_j+\Omega_j,
\label{step_cf}
\end{equation}
where
$
\Omega_j=\Omega(\varepsilon,x_j,a(x_j)).
$
By property {\rm P3)},
$$
\|\Omega_j\|\le\varpi\varepsilon^{1+\eta}\|a(x_j)\|.
$$
Since $\mathcal F^{-1}(x)$ is bounded on $\overline{B_r(0)}$, there exists $M_F>0$ such that $\|\mathcal F^{-1}(x)\|\le M_F$ for all $x\in\overline{B_r(0)}$. Hence, by~\eqref{ax} and~\eqref{gamma},
\begin{equation}
\|a(x_j)\|\le M_F\gamma(x_j)\|x_j\|=M_F\lambda\|x_j\|\left|\ln\frac{\|x_j\|}{r}\right|^\alpha.
\label{a_est}
\end{equation}
Moreover,
$
\|x\|\left|\ln\frac{\|x\|}{r}\right|^\alpha\le r\left(\frac{\alpha}{e}\right)^\alpha\,\text{for all }0<\|x\|<r,
$
so that $\|a(x_j)\|$ is uniformly bounded in $B_r(0)$. By property {\rm P3)} and~\eqref{a_est},
$
\|\Omega_j\|\le\varpi M_F\varepsilon^{1+\eta}\gamma(x_j)\|x_j\|.
$
Thus, for any $\theta\in(0,1)$, choosing $\varepsilon>0$ sufficiently small so that
$
\varpi M_F\varepsilon^\eta\le\theta,
$
we obtain
\begin{equation}
\|\Omega_j\|\le\theta\varepsilon\gamma(x_j)\|x_j\|.
\label{Omega_relative}
\end{equation}
It follows from~\eqref{step_cf} and~\eqref{Omega_relative} that
$$
\|x_{j+1}\|\le\left(\left|1-\varepsilon\gamma(x_j)\right|+\theta\varepsilon\gamma(x_j)\right)\|x_j\|.
$$
If
$
0<\varepsilon\gamma(x_j)\le1,
$
then
\begin{equation}
\begin{aligned}
\|x_{j+1}\|&\le\left(1-(1-\theta)\varepsilon\gamma(x_j)\right)\|x_j\|\\
&=\left(1-(1-\theta)\varepsilon\lambda\left|\ln\frac{\|x_j\|}{r}\right|^\alpha\right)\|x_j\|.
\end{aligned}
\label{x_step_contr}
\end{equation}
In particular, $\|x_{j+1}\|<\|x_j\|$ whenever $0<\varepsilon\gamma(x_j)\le1$.

Define
\begin{equation}
\delta_\varepsilon=re^{-\left(\frac{1}{\varepsilon\lambda}\right)^{1/\alpha}},
\label{delta_eps}
\end{equation}
which implies that if $\delta_\varepsilon\le\|x_j\|<r$, then $0<\varepsilon\gamma(x_j)\le1$, and~\eqref{x_step_contr} yields
$$
\|x_{j+1}\|\le\Big(1-(1-\theta)\varepsilon\lambda\Big|\ln\frac{\|x_j\|}{r}\Big|^\alpha\Big)\|x_j\|<\|x_j\|.
$$

For $0<\|x_j\|\le\delta_\varepsilon$, let
$
q=\ln\frac{\delta_\varepsilon}{\|x_j\|}\ge0,
$
so that $\|x_j\|=\delta_\varepsilon e^{-q}$. 

From~\eqref{delta_eps},
$
\left|\ln\frac{\|x_j\|}{r}\right|=\left(\frac{1}{\varepsilon\lambda}\right)^{1/\alpha}+q,
$
and therefore
$$
\varepsilon\gamma(x_j)=\left(1+(\varepsilon\lambda)^{1/\alpha}q\right)^\alpha\le1+\alpha(\varepsilon\lambda)^{1/\alpha}q,
$$
where we have used $(1+s)^\alpha\le1+\alpha s$ for $s\ge0$ and $\alpha\in(0,1]$. Since $\varepsilon\gamma(x_j)\ge1$, it follows from~\eqref{Omega_relative} that
$$
\begin{aligned}
\|x_{j+1}\|&\le\left((1+\theta)\varepsilon\gamma(x_j)-1\right)\|x_j\|\\
&\le\delta_\varepsilon e^{-q}\left(\theta+(1+\theta)\alpha(\varepsilon\lambda)^{1/\alpha}q\right).
\end{aligned}
$$
Taking $\varepsilon>0$ sufficiently small so that
$
(1+\theta)\alpha(\varepsilon\lambda)^{1/\alpha}\le1,
$
we obtain
$$
\|x_{j+1}\|\le\delta_\varepsilon e^{-q}(\theta+q)\le\delta_\varepsilon e^{-q}(1+q)\le\delta_\varepsilon.
$$
The case $x_j=0$ is immediate since $a(0)=0$. Thus, $\overline{B_{\delta_\varepsilon}(0)}$ is invariant at the sampling instants.

To estimate the decay rate before the trajectory enters the invariant neighborhood $B_{\delta_\varepsilon}(0)$ at the sampling instants, we introduce
$
y_j=\left|\ln\frac{\|x_j\|}{r}\right|.
$
For any $\|x_j\|\in(\delta_\varepsilon,r)$, estimate~\eqref{x_step_contr} gives
$$
y_{j+1}\ge y_j-\ln\left(1-(1-\theta)\varepsilon\lambda y_j^\alpha\right).
$$
Since $\varepsilon\lambda y_j^\alpha\le1$ and $\theta\in(0,1)$, we have $(1-\theta)\varepsilon\lambda y_j^\alpha\in[0,1)$. Using $-\ln(1-s)\ge s$ for $s\in[0,1)$, we obtain
\begin{equation}
y_{j+1}\ge y_j+(1-\theta)\varepsilon\lambda y_j^\alpha.
\label{y_rec}
\end{equation}

Define $
\bar\lambda=\frac{(1-\theta)\lambda}{2^\alpha}.
$
Let $\alpha\in(0,1)$. For sufficiently small $\varepsilon$ such that
$
\varepsilon\lambda(1-\theta)^\alpha\le1,
$
we have
$
\left|\ln\frac{\delta_\varepsilon}{r}\right|^{1-\alpha}
=\left(\frac{1}{\varepsilon\lambda}\right)^{\frac{1-\alpha}{\alpha}}
\ge(1-\theta)\lambda\varepsilon.
$
Thus, there exists a  $j_0$ such that
$
y_{j_0}^{1-\alpha}\ge(1-\theta)\lambda\varepsilon.
$
Let $j_0$ be the first such index.

For $j\ge j_0$, as long as the sampled trajectory remains outside $B_{\delta_\varepsilon}(0)$, we have
$
\frac{(1-\theta)\varepsilon\lambda}{y_j^{1-\alpha}}\le1.
$

From~\eqref{y_rec},
$$
y_{j+1}^{1-\alpha}\ge y_j^{1-\alpha}\Big(1+\frac{(1-\theta)\varepsilon\lambda}{y_j^{1-\alpha}}\Big)^{1-\alpha}.
$$
Using
$
(1+s)^{1-\alpha}-1\ge\frac{(1-\alpha)s}{(1+s)^\alpha}
$
for $s\ge0$, we obtain
$$
y_{j+1}^{1-\alpha}-y_j^{1-\alpha}\ge(1-\alpha)\bar\lambda\varepsilon.
$$
Consequently,
$
y_j\ge\Big(y_{j_0}^{1-\alpha}+(1-\alpha)\bar\lambda(t_j-t_{j_0})\Big)^{\frac{1}{1-\alpha}},
$
as long as the sampled trajectory remains outside $B_{\delta_\varepsilon}(0)$. Therefore,
\begin{equation}
\|x_j\|\le r e^{-\Big(y_{j_0}^{1-\alpha}+(1-\alpha)\bar\lambda(t_j-t_{j_0})\Big)^{\frac{1}{1-\alpha}}}.
\label{x_rate_alpha}
\end{equation}


If $\alpha=1$, estimate~\eqref{y_rec} yields
$
y_{j+1}\ge\left(1+(1-\theta)\varepsilon\lambda\right)y_j.
$
For sufficiently small $\varepsilon$ such that $(1-\theta)\varepsilon\lambda\le1$, we have
$
\frac{\ln\left(1+(1-\theta)\varepsilon\lambda\right)}{\varepsilon}
\ge\frac{(1-\theta)\lambda}{2}
=\bar\lambda.
$
Hence,
$
y_j\ge y_0e^{\bar\lambda t_j},
$
and therefore
\begin{equation}
\|x_j\|\le re^{-y_0e^{\bar\lambda t_j}},
\label{x_rate_alpha1}
\end{equation}
as long as the sampled trajectory remains outside $B_{\delta_\varepsilon}(0)$.

In both cases, these estimates imply that the sampled trajectory enters $B_{\delta_\varepsilon}(0)$ after a finite number of sampling steps. Let
$
j_\varepsilon=\min\left\{j\ge0:\|x_j\|\le\delta_\varepsilon\right\},\, T_\varepsilon=j_\varepsilon\varepsilon.
$

It remains to analyze the behavior of the $\pi_\varepsilon$-solutions of the closed-loop system for arbitrary $t\in[t_j,t_{j+1}]$. 
By property {\rm P1)}, there exists $c_u>0$ such that
$$
\sum_{k=1}^m|u_k^\varepsilon(t,x_j)|\le c_u\sum_{\ell=1}^\nu\varepsilon^{\frac{1}{\ell}-1}\|a(x_j)\|^{1/\ell}.
$$
Since
$
x(t)-x_j=\int_{t_j}^{t}\sum_{k=1}^m u_k^\varepsilon(s,x_j)f_k(x(s))\,ds,
$
and the vector fields are bounded on $K$, there exists $M_f>0$ such that $\|f_k(x)\|\le M_f$ for all $x\in K$ and $k=1,\ldots,m$. Hence, for all $t\in[t_j,t_{j+1}]$,
\begin{equation}
\|x(t)-x_j\|\le M_fc_u\sum_{\ell=1}^\nu\varepsilon^{1/\ell}\|a(x_j)\|^{1/\ell}.
\label{intersample_est}
\end{equation}

Suppose now that $\|x_j\|\le\delta_\varepsilon$. For sufficiently small $\varepsilon$, the function $s\left|\ln(s/r)\right|^\alpha$ is increasing on $(0,\delta_\varepsilon]$. Therefore, by~\eqref{a_est} and~\eqref{delta_eps},
$$
\|a(x_j)\|\le M_F\lambda\delta_\varepsilon\left|\ln\frac{\delta_\varepsilon}{r}\right|^\alpha=\frac{M_F\delta_\varepsilon}{\varepsilon}.
$$
Substituting this estimate into~\eqref{intersample_est}, we obtain
$$
\|x(t)-x_j\|\le M_fc_u\sum_{\ell=1}^\nu M_F^{1/\ell}\delta_\varepsilon^{1/\ell}.
$$
Thus, for all $t\in[t_j,t_{j+1}]$,
$$
\|x(t)\|\le\bar\delta_\varepsilon:=\delta_\varepsilon+M_fc_u\sum_{\ell=1}^\nu M_F^{1/\ell}\delta_\varepsilon^{1/\ell}.
$$
Since $\overline{B_{\delta_\varepsilon}(0)}$ is invariant at the sampling instants, once $x_j$ enters $\overline{B_{\delta_\varepsilon}(0)}$, the corresponding $\pi_\varepsilon$-solution satisfies
$$
\|x(t)\|\le\bar\delta_\varepsilon\,\text{for all }t\ge T_\varepsilon.
$$
Moreover,
$
\bar\delta_\varepsilon\to0
$
as $\varepsilon\to0$.

As long as $\delta_\varepsilon\le\|x_j\|<r$, we have $\varepsilon\gamma(x_j)\le1$, and hence
$$
\|a(x_j)\|\le M_F\gamma(x_j)\|x_j\|\le\frac{M_F}{\varepsilon}\|x_j\|.
$$
Thus,~\eqref{intersample_est} yields
$$
\|x(t)\|\le\|x_j\|+M_fc_u\sum_{\ell=1}^\nu M_F^{1/\ell}\|x_j\|^{1/\ell}
$$
for all $t\in[t_j,t_{j+1}]$. Since $\|x_j\|<r$,
\begin{equation}
\|x(t)\|\le C\|x_j\|^{1/\nu},
\label{intersample_xj}
\end{equation}
where
$
C=r^{1-\frac{1}{\nu}}+M_fc_u\sum_{\ell=1}^\nu M_F^{1/\ell}r^{\frac{1}{\ell}-\frac{1}{\nu}}.
$

For an arbitrary $t\in[0,T_\varepsilon)$, let
$
t_\varepsilon=\left[\frac{t}{\varepsilon}\right]\varepsilon,
$
where $[\cdot]$ denotes the integer part. Then
$
0\le t-t_\varepsilon<\varepsilon.
$

For $\alpha\in(0,1)$, combining~\eqref{intersample_xj} and~\eqref{x_rate_alpha}, we obtain
$$
\|x(t)\|\le Cr^{1/\nu}
e^{-\frac{1}{\nu}\Big(y_{j_0}^{1-\alpha}+(1-\alpha)\bar\lambda(t_\varepsilon-t_{j_0})\Big)^{\frac{1}{1-\alpha}}},
$$
for $t_\varepsilon\ge t_{j_0}$. Since $0\le t-t_\varepsilon<\varepsilon$ and
$
(t_\varepsilon-t_{j_0})_+^{\frac{1}{1-\alpha}}
\ge
2^{-\frac{\alpha}{1-\alpha}}t^{\frac{1}{1-\alpha}}
-(t_{j_0}+\varepsilon)^{\frac{1}{1-\alpha}},
$
we obtain
$$
\|x(t)\|\le C_\varepsilon e^{-\sigma t^{\frac{1}{1-\alpha}}},
\, t\in[0,T_\varepsilon),
$$
where
$
\sigma=\frac{2^{-\frac{\alpha}{1-\alpha}}}{\nu}
\left((1-\alpha)\bar\lambda\right)^{\frac{1}{1-\alpha}},
$
and
$
C_\varepsilon=Cr^{1/\nu}
e^{\frac{1}{\nu}\left((1-\alpha)\bar\lambda(t_{j_0}+\varepsilon)\right)^{\frac{1}{1-\alpha}}}.
$
Thus, the $\pi_\varepsilon$-solution exhibits power-type hyperexponential decay on $[0,T_\varepsilon)$.

If $\alpha=1$, combining~\eqref{x_rate_alpha1} and~\eqref{intersample_xj}, we obtain
$$
\|x(t)\|\le Cr^{1/\nu}e^{-\frac{y_0}{\nu}e^{\bar\lambda t_\varepsilon}}.
$$
Since $t_\varepsilon>t-\varepsilon$ and $\varepsilon<\hat\varepsilon$,
$$
\|x(t)\|\le Cr^{1/\nu}e^{-\sigma_0e^{\bar\lambda t}},
\qquad t\in[0,T_\varepsilon),
$$
where
$
\sigma_0=\frac{y_0}{\nu}e^{-\bar\lambda\hat\varepsilon}.
$
Thus, the $\pi_\varepsilon$-solution exhibits double-exponential decay on $[0,T_\varepsilon)$.

It remains to specify the admissible range of $\varepsilon$. Let $\varepsilon_0>0$ be such that properties {\rm P1)}--{\rm P3)} and the above estimates hold for all $\varepsilon\in(0,\varepsilon_0)$. Define
$$
\hat\varepsilon=\min\Big\{\varepsilon_0,\Big(\frac{\theta}{\varpi M_F}\Big)^{1/\eta},\frac{1}{\lambda((1+\theta)\alpha)^\alpha},\frac{1}{\lambda(1-\theta)^\alpha}\Big\}.
$$
Then, for every $\varepsilon\in(0,\hat\varepsilon)$, all the preceding estimates hold. This completes the proof.
\end{proof}

\section{Numerical Simulations}\label{sec:sim}
We consider the model of a car with two off-hooked trailers described in~\cite{ZG24}.
The system has the form
\begin{equation}\label{ex_trailers}
\dot x=u_1f_1(x)+u_2f_2(x),\quad x\in\mathbb R^5,\;u\in\mathbb R^2,
\end{equation}
where $x_1, x_2$ are the coordinates of the center of the front axle, $x_3$ is the orientation angle of the car,
$x_4$ and $x_5$ are the orientation angles of the first and the second trailers, respectively.
The controls $u_1$ and $u_2$ denote the driving velocity and steering velocity of the front axle, respectively.
The vector fields of~\eqref{ex_trailers} are given by
\begin{equation*}
\begin{aligned}
f_1(x)&=\big(\cos x_3,\ \sin x_3,\ 0,\ 
{\sin(x_3-x_4)}/{\ell_1},\ f_{15}(x)\big)^\top,\\ 
f_2(x)&=\big(0,\ 0,\ 1,\ 
-{d_0\cos(x_3-x_4)}/{\ell_1},\ f_{25}(x)\big)^\top,\\
f_{15}(x)&={\sin(x_3-x_5)}/{\ell_2}\\
&\quad+{(d_1+\ell_1)\sin(x_4-x_3)\cos(x_4-x_5)}/{\ell_1\ell_2},\\
f_{25}(x)&={d_0(d_1+\ell_1)\cos(x_3-x_4)\cos(x_4-x_5)}/{\ell_1\ell_2}\\
&\quad-{d_0\cos(x_3-x_5)}/{\ell_2}.
\end{aligned}
\end{equation*}
Here, $d_0$ and $d_1$ are the hook lengths, and $\ell_1$ and $\ell_2$ are the link lengths. In the simulations, we assume that all variables and parameters are dimensionless and set $d_0=d_1=0.1$ and $\ell_1=\ell_2=1$.

System~\eqref{ex_trailers} has degree of nonholonomy $\nu=4$ and satisfies Assumption~\ref{rank_nu}
with $n=5$ and $m=2$.
The sets of multi-indices are $S_1=\{1,2\}$, $S_2=\{(1,2)\}$, $S_3=\{(1,1,2)\}$, and $S_4=\{(1,1,1,2)\}$.
The Lie algebra rank condition (LARC) holds in a neighborhood of zero because the matrix
\begin{equation*}
\mathcal F(x)=\big(f_1,f_2,[f_1,f_2],[f_1,[f_1,f_2]],[f_1,[f_1,[f_1,f_2]]]\big)(x)
\end{equation*}
is nonsingular at $x=0$:
\begin{equation*}
\det\mathcal F(0)=-{(d_0+\ell_1)(d_0+\ell_2)(d_1+\ell_2)}/{\ell_1^4\ell_2^4}\neq0.
\end{equation*}

In~\cite{ZG24}, the nilpotent approximation of system~\eqref{ex_trailers} was stabilized by oscillating controls with a constant gain $\gamma>0$,
and exponential convergence was proved.
Here, we use controls with the state-dependent gain~\eqref{gamma}, which gives hyperexponential decay until a neighborhood of the origin is reached.

In the privileged coordinates $z\in\mathbb R^5$, the nilpotent approximation of~\eqref{ex_trailers} has the form~\cite{ZG24}:
\begin{equation}\label{sys_z}
\dot z=u_1g_1(z)+u_2g_2(z),
\end{equation}
where
$g_1(z)=\left(1,\ 0,\ -z_2,\ -z_3,\ c_2z_2^3-z_4\right)^\top$,\\ 
$g_2(z)=\left(0,\ 1,\ 0,\ -c_1z_2^2,\ -c_3z_2z_3\right)^\top.$\\
The parameters $c_1,c_2,c_3$ depend only on the lengths $d_0,d_1,\ell_1,\ell_2$.
For the lengths specified above, $(c_1,c_2,c_3)=(0.10495,\,0.0024,\,0.0112)$.
The matrix
$\mathcal G(z)=\big(g_1,g_2,[g_1,g_2],[g_1,[g_1,g_2]],[g_1,[g_1,[g_1,g_2]]]\big)(z)$
has $\det\mathcal G(z)=1$, so system~\eqref{sys_z} satisfies the LARC in $\mathbb R^5$.

We design controls of the form~\eqref{open_loop} as follows:
\begin{align}
u_1&=a_1+\big(2k_{12}\omega|a_{12}|\big)^{1/2}\cos(k_{12}\omega t)\nonumber\\
&\quad+\big(2k_{112}\omega\big)^{2/3}|a_{112}|^{1/3}\cos(k_{112}\omega t)\nonumber\\
&\quad+2\big(k_{1112}\omega\big)^{3/4}|a_{1112}|^{1/4}\cos(k_{1112}\omega t),\nonumber\\
u_2&=a_2+\big(2k_{12}\omega|a_{12}|\big)^{1/2}{\rm sign}(a_{12})\sin(k_{12}\omega t)\nonumber\\
&\quad-2\big(2k_{112}\omega\big)^{2/3}|a_{112}|^{1/3}{\rm sign}(a_{112})\cos(2k_{112}\omega t)\nonumber\\
&\quad-6\big(k_{1112}\omega\big)^{3/4}|a_{1112}|^{1/4}{\rm sign}(a_{1112})\sin(3k_{1112}\omega t),\label{uz}
\end{align}
where $\omega=2\pi/\varepsilon$ and $(k_{12},k_{112},k_{1112})=(8,6,5)$.
It was shown in~\cite{grushkovskaya2024design} that the dither components of these controls, represented as in~\eqref{open_loop}, satisfy properties P1)--P3) of Assumption~\ref{as_dither} for system~\eqref{sys_z}, with $x$ replaced by $z$ and ${\mathcal F}(x)$ replaced by ${\mathcal G}(z)$. In this case, the radius of the closed neighborhood $\overline{B_r(0)}$ in Assumption~\ref{as_dither} can be chosen arbitrarily.
The components of
$a(z)=(a_1,a_2,a_{12},a_{112},a_{1112})^\top(z)=-\gamma(z)\mathcal G^{-1}(z)z$ in~\eqref{ax} are given by:
\begin{equation}\label{a_coef}
\begin{aligned}
a_1(z)&=-\gamma(z) z_1,\, a_2(z)=-\gamma(z) z_2,\\
a_{12}(z)&=-\gamma(z)(z_1z_2+z_3),\\
a_{112}(z)&=-\gamma(z)(z_1z_3+z_4+c_1 z_2^3),\\
a_{1112}(z)&=-\gamma(z)\big((c_1+2c_2-c_3)z_1z_2^3+(c_1+3c_2)z_2^2z_3\\
&\qquad\qquad+z_1z_4+z_5\big).
\end{aligned}
\end{equation}

\begin{figure}[t]
\centering
\includegraphics[width=\columnwidth]{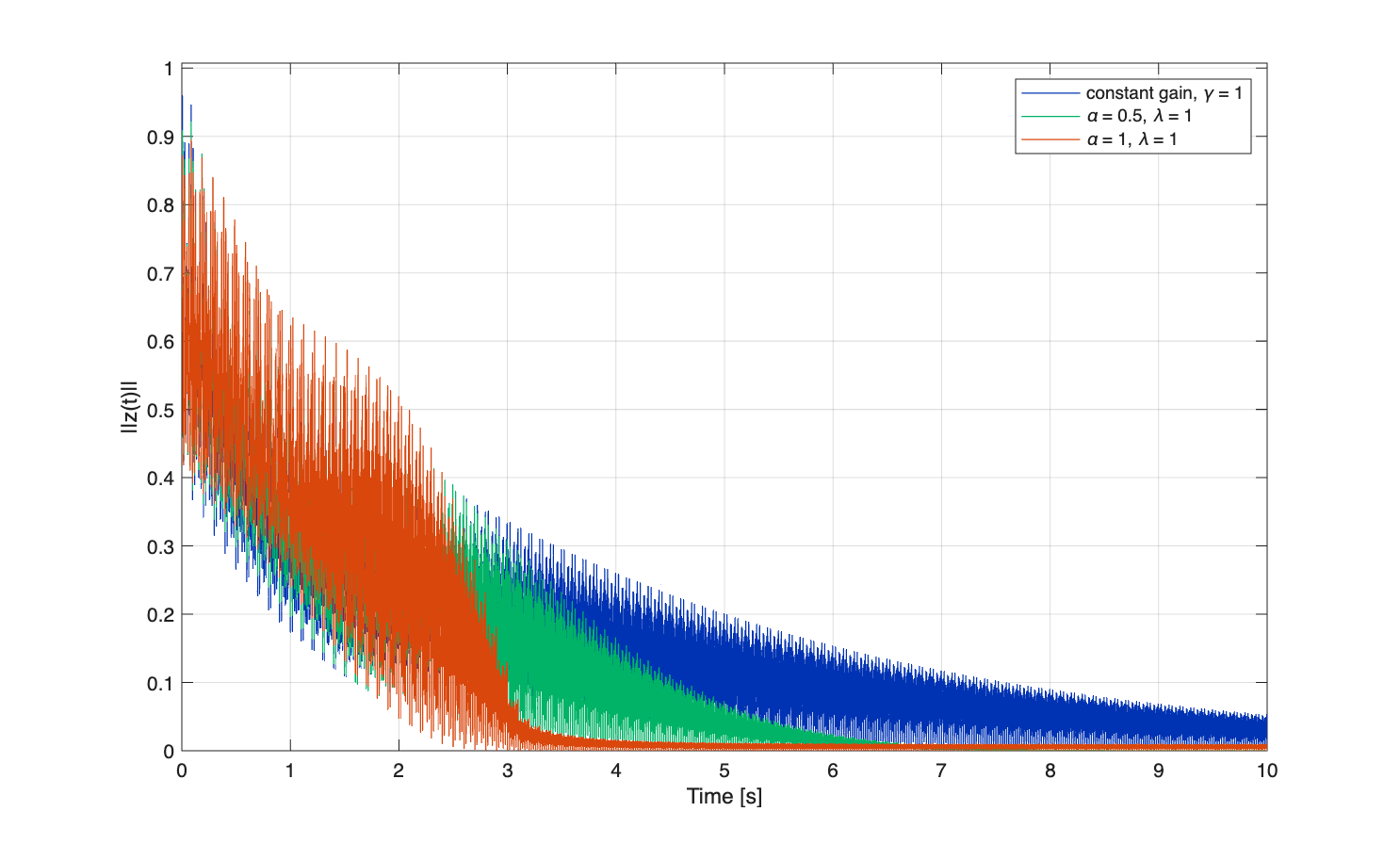}
\caption{Time plot of $\|z(t)\|$ for system~\eqref{sys_z} with gains~\eqref{eq_gains}: constant gain (blue), $\alpha=0.5$ (green), $\alpha=1$ (red).}
\label{fig_z}
\end{figure}

We perform numerical simulations of system~\eqref{sys_z} in the sense of $\pi_\varepsilon$-solutions with $\varepsilon=0.1$ using the proposed controls with three choices of gain:
\begin{equation}\label{eq_gains}
\gamma(z)\equiv1,\,
\gamma(z)=\lambda\Big|\ln\frac{\|z\|}{r}\Big|^{\alpha},\quad \alpha\in\{0.5,\,1\},
\end{equation}
where $\lambda=1$ and $r=1$.
As shown in~\cite{ZG24},
the controls with the constant gain $\gamma(z)\equiv1$ exponentially stabilize the origin of~\eqref{sys_z}, while the state-dependent gains yield the hyperexponential convergence established in Theorem~\ref{thm_main}. 
The initial state is taken as $z(0)=(0.2,-0.4,0.4,-0.2,0.2)^\top$ for all simulations in Fig.~\ref{fig_z}, so that $\|z(0)\|\approx0.663<r$.
Fig.~\ref{fig_z}illustrates that solutions $z(t)$ converge to zero faster with a state-dependent gain than with a constant gain $\gamma$, which yields exponential convergence. As follows from the comparison equation in Lemma~\ref{lem_comparison}, larger values of $\alpha\in(0,1]$ yield faster convergence, as also illustrated by comparing the simulation results for $\alpha = 1$ and $\alpha = 0.5$.

\begin{figure}[t]
\centering
\includegraphics[width=\columnwidth]{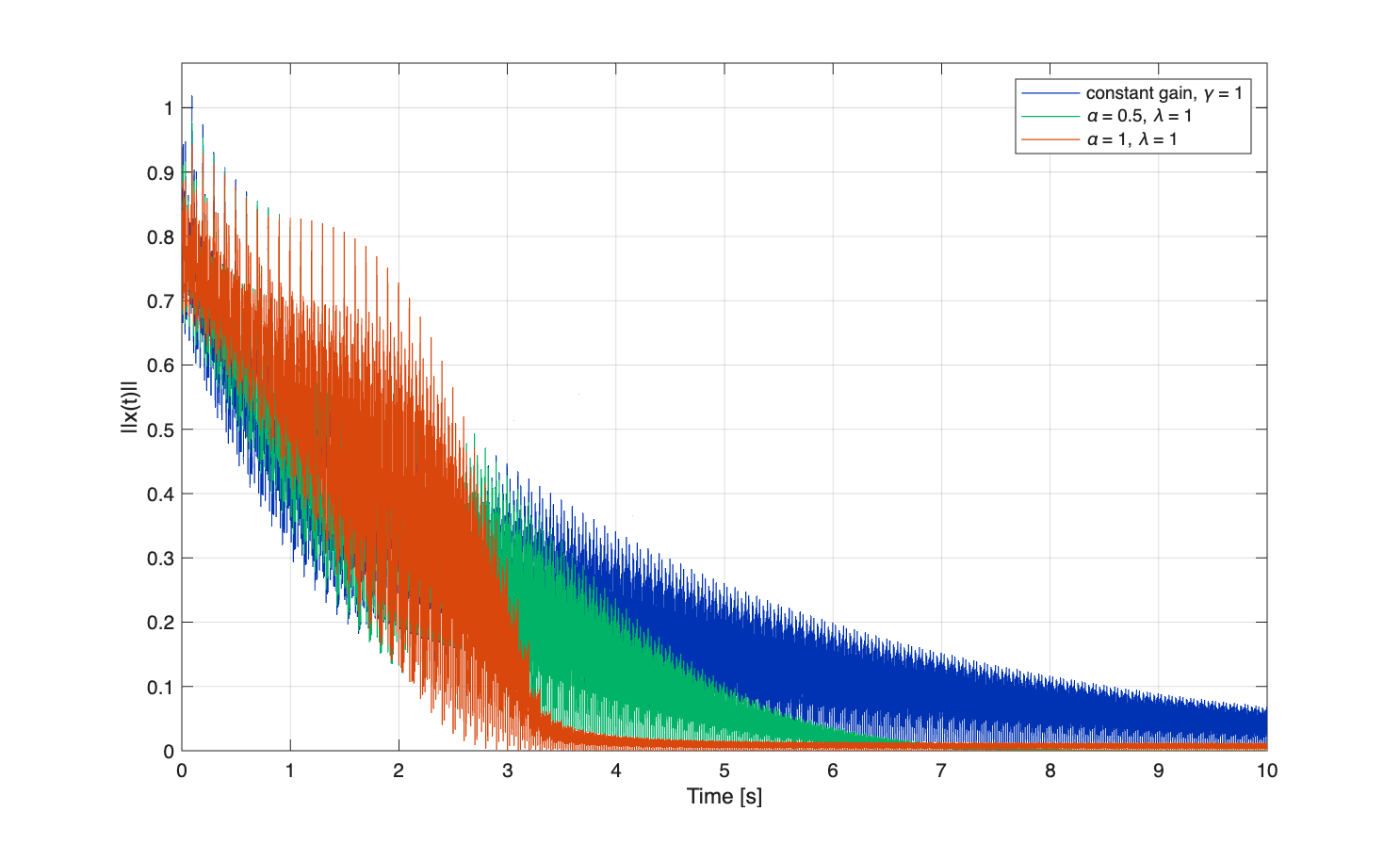}
\caption{Time plot of $\|x(t)\|$ for system~\eqref{ex_trailers} with gains~\eqref{eq_gains}: constant gain (blue), $\alpha=0.5$ (green), $\alpha=1$ (red).}
\label{fig_x}
\end{figure}

In the next simulation, the same control design is applied to the original system~\eqref{ex_trailers}.
At each sampling instant $t_j$, the state $x(t_j)$ is transformed into the privileged coordinates $z=T(x(t_j))$,
and the controls~\eqref{uz}--\eqref{eq_gains} are computed with this $z$, so that the gain depends on $\|T(x)\|$.
The resulting controls $u_1$ and $u_2$ are then applied to system~\eqref{ex_trailers} by replacing $T$ with its polynomial approximation $\hat T$ from~\cite{ZG24}, and the $\pi_\varepsilon$-solution of system~\eqref{ex_trailers} is thus computed on each interval of length $\varepsilon$. The map $\hat T$ is polynomial and has the following coordinate representation~\cite{ZG24}:
\begin{equation*}
z_1=x_1,\; z_2=x_3,\; z_3=-x_2,\; z_4=y_4,\; z_5=y_5-\beta x_3^3,
\end{equation*}
where the parameter $\beta$ depends only on the lengths $d_0,d_1,\ell_1,\ell_2$.
The initial state $x(0)$ is chosen so that $T(x(0))=z(0)$ from the previous simulation, hence $\|T(x(0))\|<r$.
Fig.~\ref{fig_x} shows the time evolution of $\|x(t)\|$ for system~\eqref{ex_trailers} under the controls described above.
We observe in Fig.~\ref{fig_x} that applying the simplified controllers designed for the nilpotentized system~\eqref{sys_z}, together with the principal part $\hat T$ of the transformation $T$, also stabilizes the original system~\eqref{ex_trailers}.


\section{Conclusion and Future Perspectives}
The transient behavior shown in Figs.~\ref{fig_z}-\ref{fig_x} illustrates that the proposed hyperexponentially stabilizing controllers produce substantially different responses from those of previously known exponentially stabilizing controllers.
Note that the logarithmic comparison equation~\eqref{comparison_log} represents a particular class of scalar dynamics generating hyperexponential convergence. We expect that this class can be generalized to generate a broader range of prescribed convergence rates. Such extensions will be addressed in future work.

\bibliographystyle{IEEEtran}
\bibliography{biblio_nonh}

\end{document}